\documentclass{ifacconf}
\usepackage{amsmath}
\usepackage{amsfonts}
\usepackage{graphicx}      
\usepackage{natbib}        
\usepackage{xcolor}
\usepackage{framed}
\usepackage{multirow}
\usepackage{lineno}
\usepackage[utf8]{inputenc}
\usepackage{textcomp}
\usepackage{booktabs}
\makeatletter
\def\endfigure{\end@float}
\def\endtable{\end@float}
\makeatother

\let\ifacconfcaptionwidth\captionwidth
\usepackage[caption=false]{subfig}
\let\captionwidth\ifacconfcaptionwidth
\usepackage{soul}
\usepackage{color}
\usepackage{comment}
\usepackage{esvect}
\usepackage[normalem]{ulem}
\usepackage{silence}
\begin{document}
\begin{frontmatter}

\title{Novel Tour Construction Heuristic for Pick-Up and Delivery Routing Problems} 
\author[First]{M. Goutham} 
\author[First]{S. Stockar} 

\address[First]{Department of Mechanical and Aerospace Engineering, The Ohio State University, 
   Columbus, OH 43210 USA (e-mail: goutham.1@osu.edu).}

\begin{abstract}                
In logistic applications that require the pickup and delivery of items, route optimization problems can be modeled as precedence constrained traveling salesperson problems. 
The combinatorial nature of this problem restricts the application of exact algorithms to small instances, and heuristics are largely preferred for tractability.
However, due to precedence constraints that restrict the order in which locations can be visited, heuristics outside of the nearest neighbor algorithm have been neglected in literature.
While the convex hull cheapest insertion heuristic is known to produce good solutions in the absence of precedence constraints, i.e., when locations can be visited in any order, it has not been adapted for pick-up and delivery considerations.
This paper presents an adapted convex hull cheapest insertion heuristic that accounts for precedence constraints and compares its solutions with the nearest neighbor heuristic using the TSPLIB benchmark data set.
The proposed algorithm is particularly suited to cases where pickups are located in the periphery and deliveries are centrally located, outperforming the Nearest Neighbor algorithm in every examined instance.

\end{abstract}


\end{frontmatter}

\section{Introduction}
Given a set of locations, the Precedence Constrained Traveling Salesperson Problem (TSP-PC) finds the shortest possible tour that respects precedence relations between locations, visiting each location exactly once and returning to the starting location \citep{mingozzi1997dynamic}.
These problems arise in a variety of real-world applications such as tool path optimization, switching energy
minimization of program compilers, automotive paint shops, and flexible manufacturing systems \citep{kucukoglu2019application, shobaki2015exact, spieckermann2004sequential}, where commodities have to be efficiently transported between pickup locations to delivery locations, and a tour is only feasible when pick up locations are visited before their respective deliveries.





Exact methods to solve TSP-PCs are limited to problems with a relatively small number of locations, typically below 100, due to their $\mathcal{NP}$-hard nature \citep{jamal2017solving, shobaki2015exact,salii2019revisiting}. 
Conversely, heuristics follow tour construction rules to greedily select locations, defining their order of visitation.
While tour construction heuristics generate feasible solutions quickly, they do not guarantee optimality and are generally sensitive to problem parameters.
For this reason, heuristics are typically used when solutions have to be found instantaneously \citep{xiang2008study, wong2014dynamic, markovic2015optimizing}, or to initialize exact methods for faster convergence to the optimal solution \citep{braekers2014exact}.

Due to the presence of precedence constraints, sophisticated heuristics tailored for the TSP-PC have been neglected in literature \citep{glover2001construction,taillard2022linearithmic}, and instead, the simple Nearest Neighbor (NN) greedy heuristic is commonly used \citep{bai2020efficient, edelkamp2018integrating}.
When precedence constraints do not exist, the Convex Hull Cheapest Insertion (CHCI) heuristic has been shown to produce superior solutions when compared to the NN heuristic \citep{ivanova2021methods,warburton1993worst}.
The CHCI heuristic is initiated by a subtour created from the convex hull of points, and its interior points are then progressively added to the subtour in increasing order of insertion cost ratios until the complete tour is obtained.
However, the CHCI heuristic has not been adapted to the TSP-PC because of the challenges in accounting for precedence constraints when inserting points to a subtour.

The main contribution of this paper is the extension of the CHCI algorithm to the TSP-PC
by initiating a subtour using the convex hull boundary points of only the non-delivery locations.
After reordering the subtour to originate at the depot, the subtour is then explored in both clockwise and counterclockwise directions because precedence constraints introduce a sense of direction of travel.
For every ensuing insertion, the adapted CHCI (ACHCI) heuristic maintains precedence constraints by only permitting insertions into feasible subtour segments.
This implies that pickup locations can be added anywhere in the subtour, while delivery locations can only be added in the subtour segment that has already visited its corresponding pickup location. 
Test cases that replicate real-world pickup and delivery considerations are created by adding precedence constraints to the TSPLIB benchmark instances \citep{reinhelt2014tsplib}.
After comparing the performance of the ACHCI and NN heuristics, their solutions are used to warm start an exact solver to study convergence.

\section{Precedence Constrained TSP}\label{NE-TSP}
\begin{figure}[b]
    \centering
        \includegraphics[trim =0mm 0mm 0mm 0mm, clip, width=0.8\linewidth]{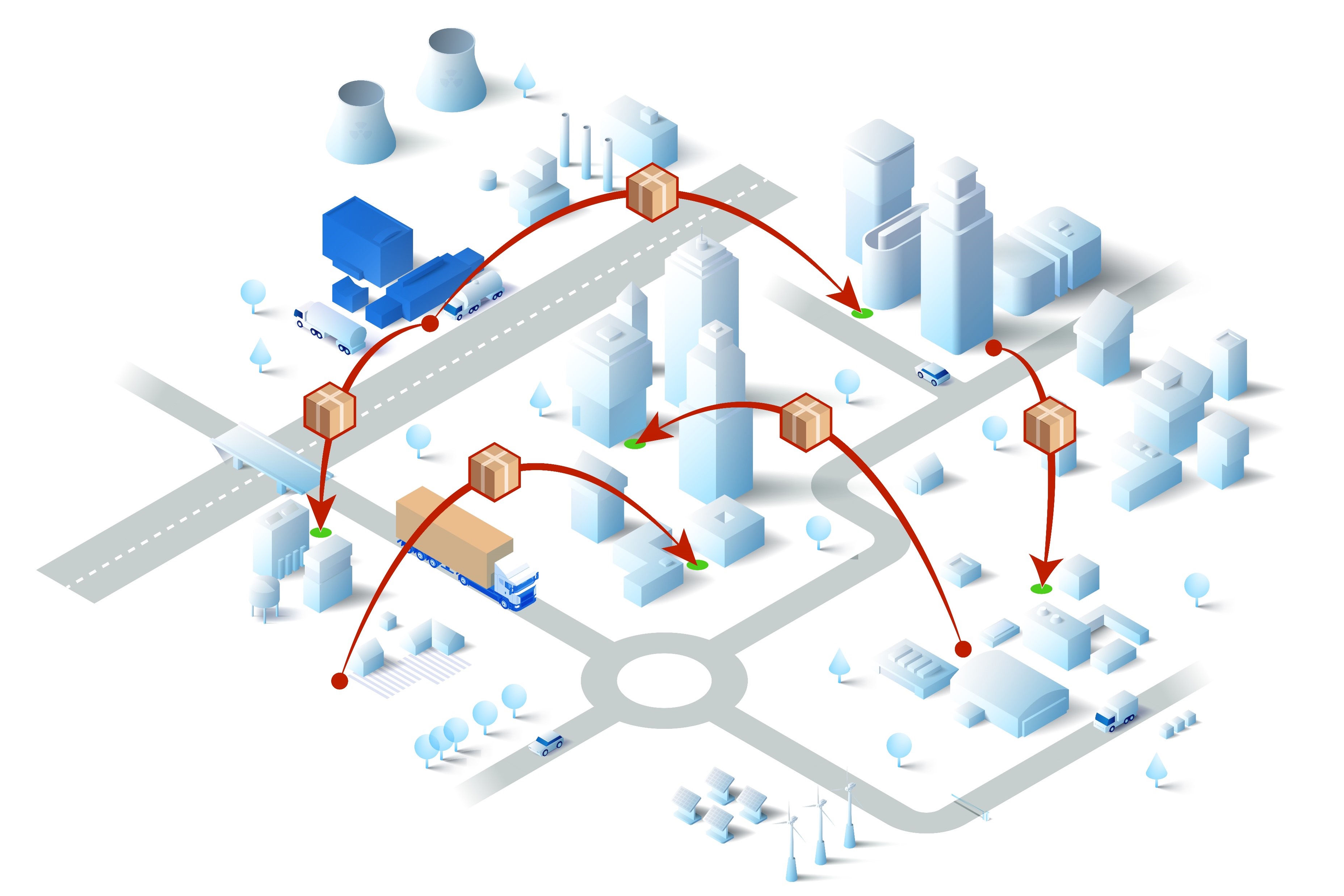}
  \caption{TSP-PC as applied to a delivery truck}
  \label{1: delivery truck} 
\end{figure}

Consider $n$ material handling tasks to be completed by an agent, with the different commodities represented by the set $\mathcal{H}:=\{h_1,h_2,...,h_n\}$.
The set of paired parent and child locations are defined by $\mathcal{V}^P:=\{1,2,...,n\}$ and $\mathcal{V}^D:=\{n+1,n+2,...,2n\}$, respectively.
A commodity $m\in\mathcal{H}$ picked up at $i\in\mathcal{V}^P$ is associated with payload $q_{im}$ and is paired with a delivery location $i+n \in \mathcal{V}^D$, such that $q_{im}+q_{i+n,m}=0$.
The start and end locations of the agent are at the depot, defined by the identically located nodes $\{0,2n+1\}$.
Define $\mathcal{V}:=\mathcal{V}^P\cup \mathcal{V}^D$,  and $\overline{\mathcal{V}}:=\mathcal{V}\cup\{0,2n+1\}$ to define the graph representation as $\mathcal{G}:=(\overline{\mathcal{V}},\mathcal{E})$, where $\mathcal{E}:=\{(i,j)\in \overline{\mathcal{V}}\times\overline{\mathcal{V}}: i\neq j\}$ denotes the set of edges.
The TSP-PC is formulated in Eq. (\ref{eq:ProbForm}) below:
\begin{subequations}  \label{eq:ProbForm}
\allowdisplaybreaks
\begin{align}  
    \label{eq:obj}      &J  = \min_{ \substack{x_{ij}\\ }} ~\sum_{(ij) \in \mathcal{E} } C_{ij} x_{ij}\\
    \label{eq:binary}   \textrm{s.t.} \quad &x_{ij}\in \{0,1\}\quad \forall (i,j)\in \mathcal{E}\\
    \label{eq:depot+}   &\sum_{j\in \mathcal{V}^p} x_{0j} = 1   \\
    \label{eq:depot-}   &\sum_{i\in \mathcal{V}^D} x_{i,2n+1}  = 1     \\ 
    \label{eq:once}     &\sum_{(ij) \in \mathcal{E}} x_{ij}  = 1           \\ 
    \label{eq:through}  &\sum_{i\in \mathcal{V}} x_{ij}= \sum_{k \in \mathcal{V}} x_{jk}        \quad\forall j \in \mathcal{V}\\
    \label{eq:subtour}  &\sum_{i,j\in \mathcal{S}} x_{ij} \leq |\mathcal{S}|-1 ~~ \forall \mathcal{S} \subset \overline{\mathcal{V}} : 2\leq |\mathcal{S}| \leq |\overline{\mathcal{V}}| -1\\
    \label{eq:cargo init}   &y_{0m} = 0          \quad\forall m\in\mathcal{H}\\
    \label{eq:cargoEvolution}
                            &\sum_{j\in \mathcal{V}}x_{ij}y_{jm} = y_{im} + \sum_{j\in \mathcal{V}}x_{ij}q_{jm}    ~\forall m\in\mathcal{H}, i \in \mathcal{V}\\
    \label{eq:cargoPrecedence}   &y_{im} \geq 0   \quad\forall m\in\mathcal{H}, i \in \mathcal{V}
        \end{align} 
\end{subequations}
The cost to travel between each node pair $(i,j)\in \overline{\mathcal{V}}$ is defined as parameter $ C_{ij}\in \mathbb{R}^+$ which is to be minimized over the tour, as captured by Eq. (\ref{eq:obj}).
Binary variables $x_{ij}$ are used to indicate whether the agent uses edge $(i,j)\in \mathcal{E}$ in the tour.
The agent must start and end at the depot, as specified by Eq. (\ref{eq:depot+}) and (\ref{eq:depot-}) respectively. 
Additionally, the agent is permitted to visit each location only once, as enforced by Eq. (\ref{eq:once}), and must leave the location after completing the visit, as defined by Eq. (\ref{eq:through}).
The subtour elimination constraint of Eq. (\ref{eq:subtour}) ensures that the resulting tour visits every location in $\overline{\mathcal{V}}$.

Payload variables $y_{im}$ define the mass of commodity $m\in\mathcal{H}$ being carried by the agent as it leaves node $i\in \overline{\mathcal{V}}$.
The agent starts its tour with no payload at the depot, as specified by Eq. (\ref{eq:cargo init}) and the evolution of payload as the agent visits pickup and delivery locations is defined by Eq. (\ref{eq:cargoEvolution})
Precedence constraints for each commodity are enforced by Eq. (\ref{eq:cargoPrecedence}) which prevents the visit of a delivery location before the corresponding item has been picked up.

\begin{figure}[t]
    \centering
        \includegraphics[trim =0mm 0mm 0mm 0mm, clip, width=0.8\linewidth]{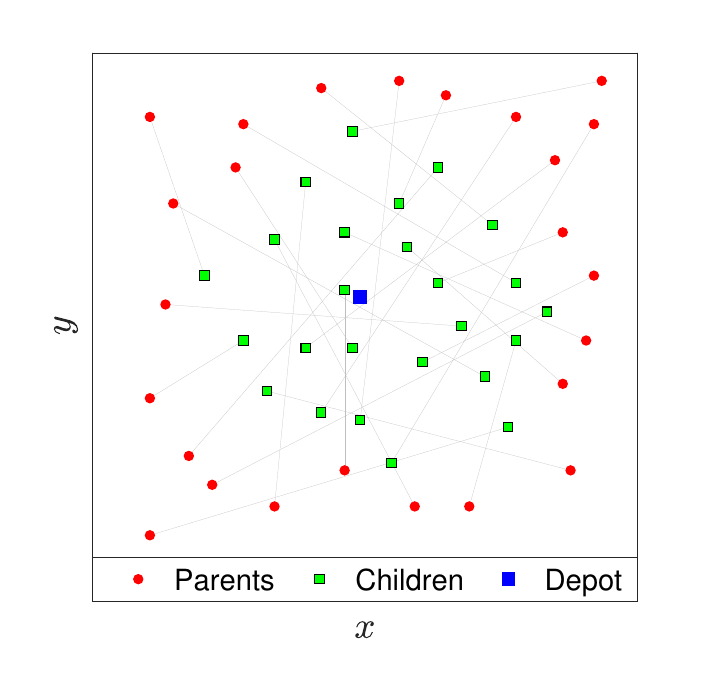}
        \vspace{-2mm}
  \caption{Generalized TSP-PC using the TSPLIB `eil51' instance with added precedence constraints}
  \label{Im3a: PC} 
\end{figure}

The problem of optimally visiting package pickup and delivery locations is shown in Fig. \ref{1: delivery truck}, where, depending on the application, $C_{ij}$ may capture distance, operational, fuel or any expense to be minimized between two locations $i,j\in \mathcal{V}$.
For an intracity problem, if $C_{ij}$ denotes distance, the tour cost of a truck that starts at the depot, picks up and delivers all the packages and returns to the depot typically amounts to tens of miles.
When applied to helicopter scheduling, this may rise to hundreds of miles \citep{fiala1992precedence}, while in the context of flexible manufacturing systems, it is often measured in feet \citep{ascheuer1993cutting}.

To generalize the analysis of the proposed ACHCI heuristic for any TSP-PC application, the popular TSPLIB benchmark instances \citep{reinhelt2014tsplib} are used with added precedence constraints, and the total Euclidean distance is minimized in this paper.
TSPLIB instance `eil51' is shown in Fig. \ref{Im3a: PC}, where precedence constraints are illustrated using grey line segments that relate parent pickup positions, shown in red, with the respective children drop-off positions, marked in green. 
The agent starts and ends its tour at the depot, marked in blue.

\section{Adapted CHCI Heuristic}


The ACHCI heuristic is initiated as the ordered boundary nodes of the convex hull of non-child locations, that is, the set $\bar{V}_0:=\mathcal{V}^P \cup \{0,2n+1\}$.
This ordered set of points defines a node sequence $T_0 := [v_1, v_2, ..., v_1]$.
For any node $v_k$ not in $T_0$, and any two consecutive nodes $v_i,v_j$ in $T_0$, the insertion cost ratio of $v_k$ with respect to $v_i,v_j$ is:
\begin{equation}  \label{eq:insertionCost}
(C_{ik} + C_{kj})/C_{ij}
\end{equation}
If the identically located depot node $\{0,2n+1\}\notin T_0$, consecutive nodes $v_i^*,v_j^*$ with the lowest insertion cost ratio are first found in $T_0$, and the depot node is inserted between them.
To start and end the subtour at the depot, the updated subtour $T_0$ is reordered to $[0,v_j^*, ...,v_i^*, 2n+1]$.
Because a sense of direction exists due to the precedence constraints, denote the resulting subtour as $\overleftarrow{T}$.
The subsequent steps will be repeated in the other direction $\overrightarrow{T}$ as well, and the lower cost tour is selected. 

\begin{figure}[t]
    \centering
      \vspace{-3mm}
  \subfloat[Initial convex hull subtour\label{Im3b: cHull subtour} ]{%
       \includegraphics[trim =0mm 0mm 0mm 0mm, clip, width=0.62\linewidth]{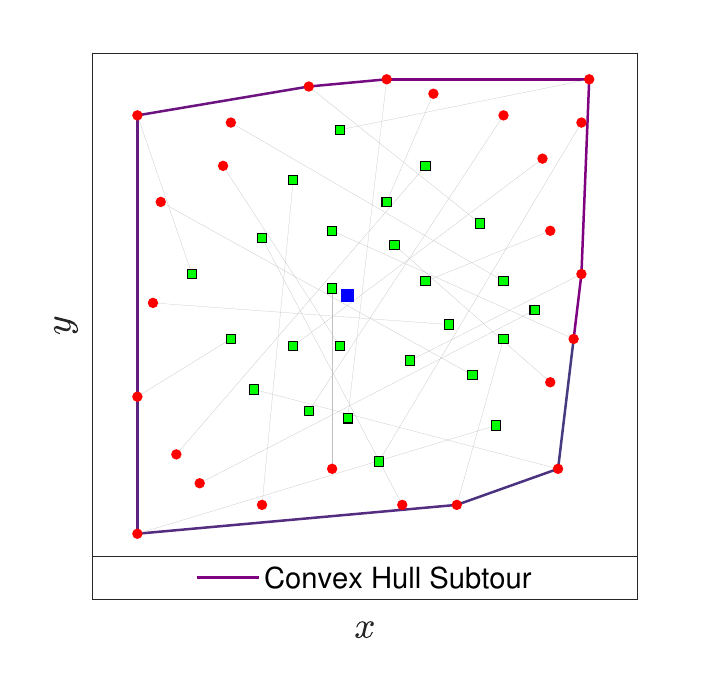}}
       \\
       \vspace{0mm}
\subfloat[Child node insertion onto the valid segment of the counterclockwise tour \label{Im4a: Insert1 CCW}]{%
       \includegraphics[trim =0mm 0mm 0mm 0mm, clip, width=0.52\linewidth]{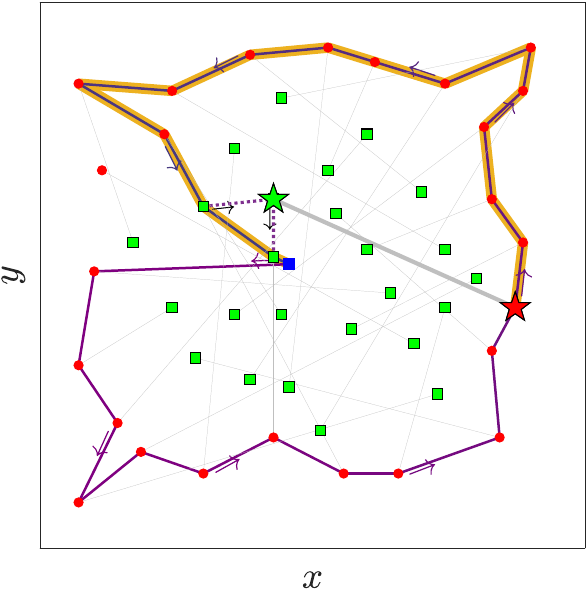}}
       \\
       \vspace{0mm}
  \subfloat[Completed counterclockwise tour \label{Im4b: Tour CCW}]{%
        \includegraphics[trim =0mm 0mm 0mm 0mm, clip, width=0.52\linewidth]{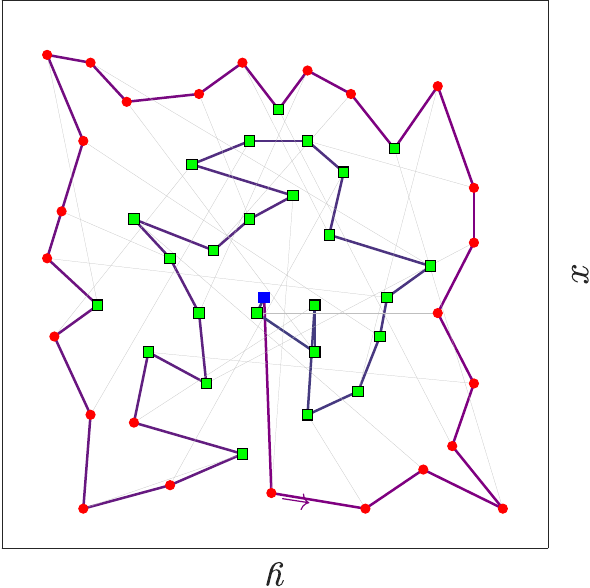}}
  \caption{Illustrative example tour of the ACHCI algorithm}
  \label{fig2} 
  \vspace{0mm}
\end{figure}

For a node $v_p\in \overleftarrow{T}$, let the segment of $\overleftarrow{T}$ that has already visited $v_p$ be denoted by $\overleftarrow{T}_{p^+} := [v_p,v_{p+1},..., 2n+1]$.
For a candidate insertion node $v_k \notin \overleftarrow{T}$, the feasible partition of $\overleftarrow{T}$ where $v_k$ can be inserted is the subtour segment that contains every parent of this candidate.
If the set of all parents of $v_k$, denoted by $\mathcal{V}_{k}^-$ is empty, then the candidate node $v_k$ can be inserted anywhere in the subtour.
Contrarily, if $|\mathcal{V}_{k}^-|>0$, then the insertion of $v_k \notin \overleftarrow{T}$ is only permitted in segment $\overleftarrow{T}^k = \bigcap_{p\in \mathcal{V}_{k}^-}\overleftarrow{T}_{p^+}$.

For every node $v_i \notin \overleftarrow{T}$ that is yet to be inserted to the subtour, the insertion arc given by consecutive nodes $(v_q,v_r) \in \overleftarrow{T}^i$ is found that minimizes insertion cost ratio given by Eq. (\ref{eq:insertionCost}).
To ensure feasibility, child nodes whose parents have not yet been visited in $\overleftarrow{T}$ are assigned an infinite insertion cost.
Next, the node $v^* \notin \overleftarrow{T}$ with the lowest insertion cost ratio is inserted at its associated insertion arc.
This increments $\overleftarrow{T}$ by one node, and these steps are repeated until every node has been inserted.

Considering the initial subtour $T_0$ was assigned some arbitrary direction that resulted in $\overleftarrow{T}$, all of the tour constructing steps are also repeated after initializing subtour $T_0$ in the opposite direction, forming another subtour $\overrightarrow{T}$.
Thus, two complete tours $\overleftarrow{T}$ and $\overrightarrow{T}$ are obtained for the TSP-PC, and the minimum cost tour is selected.


For the convex hull points shown in Fig. \ref{Im3b: cHull subtour}, the tour with the counterclockwise direction is shown in Fig. \ref{Im4a: Insert1 CCW} at an instance when a child node, marked by a green star, is inserted.
Its associated parent is marked by the red star, and the valid partition of the tour is highlighted, being the tour segment that has already visited the parent.
Notice also that at this instance, the sub-tour $\overrightarrow{T}$ has already been incremented by other parent nodes, but not all parent nodes have been inserted.
After inserting all the remaining nodes, the completed Hamiltonian tour with some cost $\overrightarrow{J}$ is shown in Fig. \ref{Im4b: Tour CCW}.
These steps are also repeated after initiating $T_0$ in the clockwise direction, to obtain a tour of some cost $\overleftarrow{J}$. 
The tour with lower cost between $\overrightarrow{J}$ \& $\overleftarrow{J}$ is selected as the ACHCI heuristic tour.

\begin{figure*}[t!]
\centering
\begin{minipage}{0.7\textwidth}
  \subfloat[Instance `eil76' with central children \label{a: cChildren}]{%
       \includegraphics[trim =0mm 0mm 0mm 0mm, clip, width=0.5\linewidth]{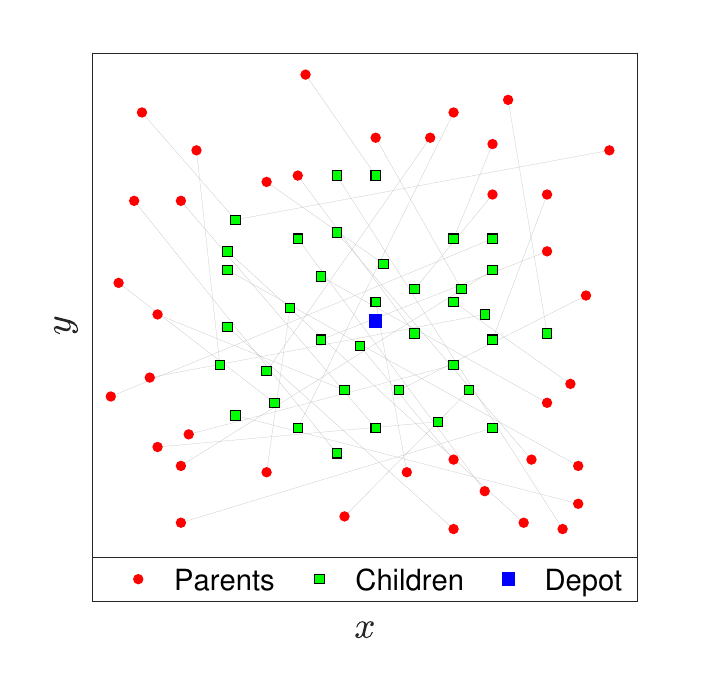}}
\hspace{10mm}
  \subfloat[ Instance `eil76' with central parents\label{b: cParents}]{%
        \includegraphics[trim =0mm 0mm 0mm 0mm, clip, width=0.5\linewidth]{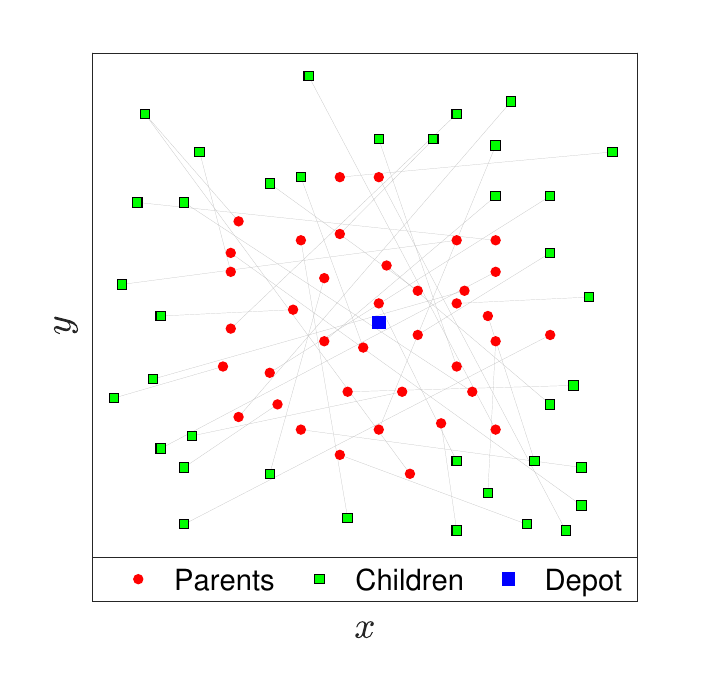}}
        \vspace{5mm}
\end{minipage}
\begin{minipage}{0.67\textwidth}
  \subfloat[Cost ratios for central children \label{d: resChildren}]{%
       \includegraphics[trim =0mm 0mm 0mm 0mm, clip, width=0.465\linewidth]{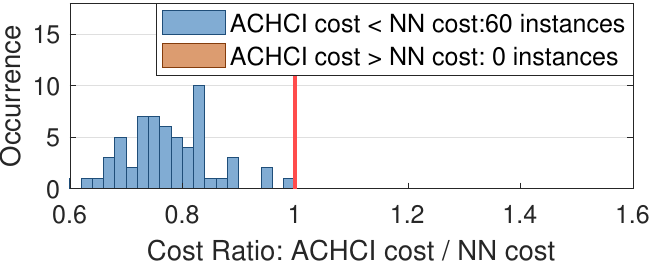}}
\hspace{16.5mm}
  \subfloat[Cost ratios for central parents\label{e: resParents  }]{%
        \includegraphics[trim =0mm 0mm 0mm 0mm, clip, width=0.465\linewidth]{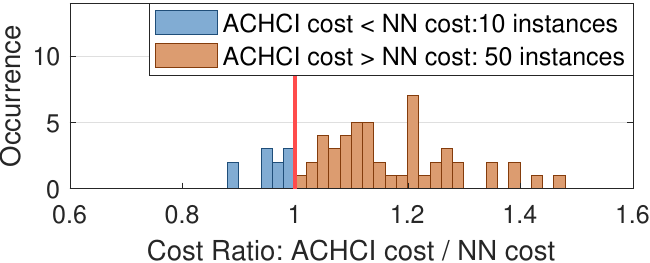}}
\end{minipage}
\captionsetup{justification=centering}
  \caption{Performance comparison of the ACHCI and NN algorithm}
  \label{PC-examples2} 
\end{figure*}

\section{Nearest Neighbor Heuristic}
The benchmark for the ACHCI algorithm is the NN which is known to produce feasible tours in real time for problems with a large number of locations \citep{charikar1997constrained}.
The NN heuristic uses a fast greedy selection rule and is commonly seen in constrained TSP literature because it can easily be modified to account for constraints in the problem formulation \citep{grigoryev2018solving,bai2020efficient}.
Starting from the depot, the NN heuristic assigns the nearest unvisited node as the next node until all nodes are contained in the tour.
By modifying this to only consider unvisited feasible nodes as the next node, it can be adapted for the TSP-PC:
\begin{enumerate}
    \item Initiate the subtour as the depot
    \item Add the closest feasible node to the end of the subtour
    \item Repeat step 2 until all the nodes are included
    \item Return to the depot
\end{enumerate}

\section{Computational Experiments}

To compare the effectiveness of the ACHCI heuristic with the NN heuristic, sufficiently diverse benchmark instances are not readily available for the TSP-PC.
For this reason, the popular TSPLIB benchmark instances \citep{reinhelt2014tsplib} are modified in a reproducible manner to create cases with precedence constraints. 

The following steps are used to create TSP-PC instances:
\begin{enumerate}
    \item Load a TSPLIB point cloud that has 2D Cartesian coordinates defined for every point. Let $n$ be the number of points.
    \item Find the centroid of the point cloud. Sort and assign indices to the points in order of increasing distance from the centroid. Let the indices be $1,2,...,n$ for the sorted points.
    \item Designate the point with index 1 as the depot.
    \item Define precedence constraints between points with pairs of indices as $(n\prec 2),(n-1\prec 3),(n-2\prec 4)$ and so on. If three nodes remain to be assigned precedence constraints, the node with the smallest index is assigned to be the child node while the other two are its parents. The generated point cloud has child nodes that are clustered closer to the centroid, with parent nodes at the periphery.
    \item Assign payload $q_{im} = 1$ at all parent locations $i\in \mathcal{V}^P$ and $q_{im} = -1$ for child nodes with only one parent, and $q_{im} = -2$ when a child node has two parents.
    \item Repeat (5) after changing direction of the precedence constraints to create the complementary case where parent nodes are closer to the centroid.
\end{enumerate}
The resulting spatial characteristic of the precedence constraints are illustrated using the modified TSPLIB `eil76' instance in Fig. \ref{a: cChildren} and \ref{b: cParents}. 

\subsection{Comparison of Heuristic Solutions}

The ACHCI and NN tours are compared for 60 TSPLIB instances with precedence constraints, with the number of points varying from 51 to 1,577.
The result of these experiments is summarized in Table \ref{table} where the first column lists the name of each TSPLIB instance, formatted with a prefix followed by a numeric value that indicates the number of points in the respective instance.
The second and third columns provide the solution costs of the NN and ACHCI heuristics along with the percentage reduction in NN cost by using the ACHCI solution.
The histogram in Fig. \ref{d: resChildren} shows that the ACHCI tour cost is lower than the NN cost in all cases with centrally located children.
However, the ACHCI heuritic does not perform as well when parent nodes are closer to the centroid, as seen in Fig. \ref{e: resParents  }.
\begin{table*}[t]
\centering
\caption{TSP-PC Tour Costs}\label{table}
\centering
    \resizebox{1.25\columnwidth}{!}{
    \begin{tabular}[t]{l|rl|rll}
    \toprule
\multicolumn{1}{c|}{TSPLIB}   & \multicolumn{2}{c|}{\textbf{Centrally Located Children}}          & \multicolumn{2}{c}{\textbf{Centrally Located Parents}}     \\
\multicolumn{1}{c|}{instance} & \multicolumn{1}{c}{NN} & \multicolumn{1}{c|}{ACHCI} & \multicolumn{1}{c}{NN} & \multicolumn{1}{c}{ACHCI}                       \\ \hline
\textit{eil51}    &5.8e+02   &4.8e+02 (-17.7\%)   &5.8e+02   &6.2e+02 (\texttt{+}\textbf{ 6.3}\%)   \\ 
\textit{t70}    &9.1e+02   &7.9e+02 (-13.3\%)   &1.1e+03   &1.2e+03 (\texttt{+}\textbf{10.6}\%)   \\ 
\textit{eil76}    &7.9e+02   &6.3e+02 (-20.4\%)   &7.3e+02   &8.3e+02 (\texttt{+}\textbf{13.3}\%)   \\ 
\textit{erlin52}    &1.1e+04   &8.3e+03 (-26.1\%)   &1.1e+04   &1.2e+04 (\texttt{+}\textbf{ 5.5}\%)   \\ 
\textit{eil101}    &9.7e+02   &7.5e+02 (-23.0\%)   &9.0e+02   &9.1e+02 (\texttt{+}\textbf{ 0.7}\%)   \\ 
\textit{rat99}    &1.8e+03   &1.5e+03 (-17.4\%)   &1.8e+03   &1.9e+03 (\texttt{+}\textbf{ 7.7}\%)   \\ 
\textit{pr76}    &2.0e+05   &1.3e+05 (-35.0\%)   &1.5e+05   &1.7e+05 (\texttt{+}\textbf{11.1}\%)   \\ 
\textit{roC100}    &3.5e+04   &2.5e+04 (-30.4\%)   &2.9e+04   &3.7e+04 (\texttt{+}\textbf{27.1}\%)   \\ 
\textit{roD100}    &3.5e+04   &2.8e+04 (-21.2\%)   &2.8e+04   &3.5e+04 (\texttt{+}\textbf{26.6}\%)   \\ 
\textit{roE100}    &3.3e+04   &2.8e+04 (-17.2\%)   &3.2e+04   &3.6e+04 (\texttt{+}\textbf{12.3}\%)   \\ 
\textit{roA100}    &3.0e+04   &2.7e+04 (-10.2\%)   &3.4e+04   &3.3e+04 (- 3.5\%)   \\ 
\textit{roB100}    &3.7e+04   &2.7e+04 (-26.1\%)   &3.3e+04   &3.8e+04 (\texttt{+}\textbf{13.8}\%)   \\ 
\textit{in105}    &2.6e+04   &1.7e+04 (-36.7\%)   &2.3e+04   &2.8e+04 (\texttt{+}\textbf{20.9}\%)   \\ 
\textit{pr107}    &7.5e+04   &7.2e+04 (- 4.9\%)   &7.2e+04   &8.6e+04 (\texttt{+}\textbf{19.0}\%)   \\ 
\textit{pr124}    &1.0e+05   &8.1e+04 (-22.1\%)   &9.5e+04   &1.1e+05 (\texttt{+}\textbf{11.3}\%)   \\ 
\textit{roB150}    &4.4e+04   &3.3e+04 (-25.8\%)   &3.9e+04   &4.0e+04 (\texttt{+}\textbf{ 1.7}\%)   \\ 
\textit{roA150}    &4.1e+04   &3.3e+04 (-21.2\%)   &4.0e+04   &3.9e+04 (- 1.5\%)   \\ 
\textit{pr136}    &1.5e+05   &1.2e+05 (-19.8\%)   &1.5e+05   &1.4e+05 (-10.7\%)   \\ 
\textit{pr144}    &9.4e+04   &7.2e+04 (-23.3\%)   &7.6e+04   &1.0e+05 (\texttt{+}\textbf{35.8}\%)   \\ 
\textit{pr152}    &1.1e+05   &8.6e+04 (-18.5\%)   &1.0e+05   &1.4e+05 (\texttt{+}\textbf{43.1}\%)   \\ 
\textit{rat195}    &3.2e+03   &2.9e+03 (-10.4\%)   &3.3e+03   &3.4e+03 (\texttt{+}\textbf{ 3.4}\%)   \\ 
\textit{bier127}    &1.8e+05   &1.4e+05 (-25.7\%)   &1.7e+05   &1.7e+05 (\texttt{+}\textbf{ 4.2}\%)   \\ 
\textit{roA200}    &5.6e+04   &3.8e+04 (-32.3\%)   &4.9e+04   &5.0e+04 (\texttt{+}\textbf{ 4.1}\%)   \\ 
\textit{roB200}    &4.9e+04   &3.6e+04 (-27.3\%)   &4.1e+04   &5.2e+04 (\texttt{+}\textbf{27.8}\%)   \\ 
\textit{rd100}    &1.5e+04   &1.0e+04 (-32.8\%)   &1.2e+04   &1.1e+04 (- 5.5\%)   \\ 
\textit{gil262}    &3.8e+03   &3.1e+03 (-19.7\%)   &3.5e+03   &4.2e+03 (\texttt{+}\textbf{21.1}\%)   \\ 
\textit{pr226}    &1.2e+05   &1.1e+05 (- 5.8\%)   &1.1e+05   &1.4e+05 (\texttt{+}\textbf{21.4}\%)   \\ 
\textit{a280}    &4.0e+03   &3.3e+03 (-17.8\%)   &4.0e+03   &3.7e+03 (- 5.8\%)   \\ 
\textit{ts225}    &2.1e+05   &1.5e+05 (-27.8\%)   &1.7e+05   &2.4e+05 (\texttt{+}\textbf{46.5}\%)   \\ 
\textit{pr264}    &7.4e+04   &7.4e+04 (- 0.5\%)   &6.8e+04   &7.6e+04 (\texttt{+}\textbf{11.2}\%)   \\ 
\textit{tsp225}    &6.0e+03   &4.3e+03 (-27.9\%)   &5.2e+03   &6.5e+03 (\texttt{+}\textbf{26.0}\%)   \\ 
\textit{pr299}    &8.1e+04   &6.1e+04 (-24.1\%)   &7.0e+04   &6.9e+04 (- 0.9\%)   \\ 
\textit{in318}    &6.2e+04   &5.2e+04 (-16.3\%)   &6.8e+04   &8.2e+04 (\texttt{+}\textbf{20.6}\%)   \\ 
\textit{in318}    &6.2e+04   &5.2e+04 (-16.3\%)   &6.8e+04   &8.2e+04 (\texttt{+}\textbf{20.6}\%)   \\ 
\textit{h130}    &9.3e+03   &7.7e+03 (-16.9\%)   &9.5e+03   &1.0e+04 (\texttt{+}\textbf{ 8.0}\%)   \\ 
\textit{u159}    &6.5e+04   &5.4e+04 (-16.8\%)   &6.0e+04   &8.4e+04 (\texttt{+}\textbf{39.2}\%)   \\ 
\textit{h150}    &9.5e+03   &8.0e+03 (-15.3\%)   &9.3e+03   &1.3e+04 (\texttt{+}\textbf{34.4}\%)   \\ 
\textit{d198}    &2.6e+04   &1.8e+04 (-30.5\%)   &2.3e+04   &2.1e+04 (-11.1\%)   \\ 
\textit{pr439}    &1.8e+05   &1.3e+05 (-27.0\%)   &1.5e+05   &1.5e+05 (- 3.4\%)   \\ 
\textit{rat575}    &1.1e+04   &8.2e+03 (-24.1\%)   &1.0e+04   &1.1e+04 (\texttt{+}\textbf{ 8.9}\%)   \\ 
\textit{rat783}    &1.5e+04   &1.0e+04 (-29.2\%)   &1.3e+04   &1.5e+04 (\texttt{+}\textbf{15.8}\%)   \\ 
\textit{rd400}    &2.4e+04   &1.8e+04 (-23.9\%)   &2.1e+04   &2.7e+04 (\texttt{+}\textbf{28.9}\%)   \\ 
\textit{fl417}    &1.9e+04   &1.7e+04 (-11.0\%)   &1.8e+04   &2.3e+04 (\texttt{+}\textbf{29.4}\%)   \\ 
\textit{pcb442}    &7.7e+04   &5.9e+04 (-22.7\%)   &8.7e+04   &8.2e+04 (- 5.0\%)   \\ 
\textit{d493}    &5.0e+04   &4.1e+04 (-18.0\%)   &4.7e+04   &5.0e+04 (\texttt{+}\textbf{ 5.1}\%)   \\ 
\textit{pr1002}    &4.1e+05   &3.1e+05 (-24.3\%)   &4.0e+05   &4.4e+05 (\texttt{+}\textbf{10.1}\%)   \\ 
\textit{u574}    &5.6e+04   &4.5e+04 (-20.6\%)   &5.9e+04   &6.9e+04 (\texttt{+}\textbf{16.9}\%)   \\ 
\textit{p654}    &6.8e+04   &5.3e+04 (-22.5\%)   &5.4e+04   &6.7e+04 (\texttt{+}\textbf{23.7}\%)   \\ 
\textit{d657}    &8.1e+04   &5.7e+04 (-28.8\%)   &7.0e+04   &8.5e+04 (\texttt{+}\textbf{20.6}\%)   \\ 
\textit{u724}    &6.9e+04   &5.2e+04 (-25.0\%)   &5.8e+04   &7.2e+04 (\texttt{+}\textbf{24.7}\%)   \\ 
\textit{u1060}    &3.6e+05   &2.7e+05 (-25.9\%)   &3.4e+05   &3.9e+05 (\texttt{+}\textbf{15.2}\%)   \\ 
\textit{vm1084}    &4.2e+05   &2.8e+05 (-31.5\%)   &3.8e+05   &4.0e+05 (\texttt{+}\textbf{ 3.6}\%)   \\ 
\textit{nrw1379}    &8.0e+04   &6.6e+04 (-17.0\%)   &8.2e+04   &8.9e+04 (\texttt{+}\textbf{ 9.4}\%)   \\ 
\textit{pcb1173}    &9.7e+04   &7.1e+04 (-26.9\%)   &8.7e+04   &8.5e+04 (- 2.0\%)   \\ 
\textit{d1291}    &7.7e+04   &6.5e+04 (-16.1\%)   &8.8e+04   &9.3e+04 (\texttt{+}\textbf{ 6.3}\%)   \\ 
\textit{rl1304}    &4.7e+05   &3.1e+05 (-33.9\%)   &4.7e+05   &5.4e+05 (\texttt{+}\textbf{13.9}\%)   \\ 
\textit{rl1323}    &5.8e+05   &3.5e+05 (-40.1\%)   &4.6e+05   &5.6e+05 (\texttt{+}\textbf{21.2}\%)   \\ 
\textit{fl1400}    &4.3e+04   &2.9e+04 (-31.7\%)   &3.1e+04   &3.5e+04 (\texttt{+}\textbf{12.6}\%)   \\ 
\textit{u1432}    &2.3e+05   &1.8e+05 (-20.1\%)   &2.2e+05   &2.4e+05 (\texttt{+}\textbf{ 9.0}\%)   \\ 
\textit{fl1577}    &3.8e+04   &2.6e+04 (-30.1\%)   &3.4e+04   &4.8e+04 (\texttt{+}\textbf{38.6}\%)   \\ 
    \bottomrule
    \end{tabular}}
\end{table*}

\begin{table*}[t]
\centering
\caption{TSP-PC Tour Costs for Instances with Centrally Located Children}\label{table1}
\begin{minipage}{1\textwidth}
\centering
    \begin{tabular}[t]{l|rl|rll}
    \toprule
\multicolumn{1}{c|}{TSPLIB}   & \multicolumn{2}{c|}{\textbf{Heuristic}}             & \multicolumn{3}{c}{\textbf{Exact}}                                                                        \\
\multicolumn{1}{c|}{instance} & \multicolumn{1}{c}{NN} & \multicolumn{1}{c|}{ACHCI} & \multicolumn{1}{c}{Standalone} & \multicolumn{1}{c}{NN warm start} & \multicolumn{1}{c}{ACHCI warm start} \\ \hline
\textit{eil51}    &5.83e+02   &4.80e+02 (-17.7\%)   &4.73e+02   &4.54e+02 (- 4.0\%)   &4.64e+02 (- 2.0\%)   \\ 
\textit{erlin52}    &1.12e+04   &8.25e+03 (-26.1\%)   &8.54e+03   &9.39e+03 (\texttt{+}\textbf{10.1}\%)   &8.14e+03 (- 4.6\%)   \\ 
\textit{t70}    &9.10e+02   &7.89e+02 (-13.3\%)   &2.20e+03   &8.29e+02 (-62.3\%)   &7.63e+02 (-65.3\%)   \\ 
\textit{eil76}    &7.93e+02   &6.32e+02 (-20.4\%)   &2.41e+03   &7.22e+02 (-70.0\%)   &6.28e+02 (-73.9\%)   \\ 
\textit{pr76}    &1.96e+05   &1.27e+05 (-35.0\%)   &4.80e+05   &1.66e+05 (-65.4\%)   &1.27e+05 (-73.5\%)   \\ 
\textit{rat99}    &1.79e+03   &1.48e+03 (-17.4\%)   &    -    &1.74e+03 (  -  ) &1.46e+03 (  -  )\\ 
\textit{roC100}    &3.52e+04   &2.45e+04 (-30.4\%)   &    -    &3.45e+04 (  -  ) &2.44e+04 (  -  )\\ 
\textit{roD100}    &3.50e+04   &2.75e+04 (-21.2\%)   &    -    &3.34e+04 (  -  ) &2.74e+04 (  -  )\\ 
\textit{roE100}    &3.34e+04   &2.77e+04 (-17.2\%)   &    -    &3.09e+04 (  -  ) &2.77e+04 (  -  )\\ 
\textit{roA100}    &3.01e+04   &2.70e+04 (-10.2\%)   &    -    &2.86e+04 (  -  ) &2.68e+04 (  -  )\\
    \bottomrule
    \end{tabular}
\end{minipage}%

\vspace{5mm}
\caption{TSP-PC Tour Costs for Instances with Centrally Located Parents}\label{table2}
\begin{minipage}{1\textwidth}
\centering
    \begin{tabular}[t]{l|rl|rll}
    \toprule
\multicolumn{1}{c|}{TSPLIB}   & \multicolumn{2}{c|}{\textbf{Heuristic}}             & \multicolumn{3}{c}{\textbf{Exact}}                                                                        \\
\multicolumn{1}{c|}{instance} & \multicolumn{1}{c}{NN} & \multicolumn{1}{c|}{ACHCI} & \multicolumn{1}{c}{Standalone} & \multicolumn{1}{c}{NN warm start} & \multicolumn{1}{c}{ACHCI warm start} \\ \hline
\textit{eil51}    &5.82e+02   &6.18e+02 (\texttt{+}\textbf{ 6.3}\%)   &1.09e+03   &5.20e+02 (-52.4\%)   &4.74e+02 (-56.6\%)   \\ 
\textit{erlin52}    &1.11e+04   &1.17e+04 (\texttt{+}\textbf{ 5.5}\%)   &9.34e+03   &9.96e+03 (\texttt{+}\textbf{ 6.6}\%)   &8.21e+03 (-12.2\%)   \\ 
\textit{t70}    &1.08e+03   &1.20e+03 (\texttt{+}\textbf{10.6}\%)   &2.57e+03   &1.03e+03 (-60.0\%)   &1.10e+03 (-57.0\%)   \\ 
\textit{eil76}    &7.31e+02   &8.29e+02 (\texttt{+}\textbf{13.3}\%)   &2.08e+03   &7.24e+02 (-65.1\%)   &7.12e+02 (-65.7\%)   \\ 
\textit{pr76}    &1.54e+05   &1.71e+05 (\texttt{+}\textbf{11.1}\%)   &4.34e+05   &1.50e+05 (-65.4\%)   &1.63e+05 (-62.4\%)   \\ 
\textit{rat99}    &1.80e+03   &1.94e+03 (\texttt{+}\textbf{ 7.7}\%)   &    -    &1.79e+03 (  -  ) &1.92e+03 (  -  )\\ 
\textit{roC100}    &2.88e+04   &3.66e+04 (\texttt{+}\textbf{27.1}\%)   &    -    &2.88e+04 (  -  ) &3.52e+04 (  -  )\\ 
\textit{roD100}    &2.80e+04   &3.55e+04 (\texttt{+}\textbf{26.6}\%)   &    -    &2.78e+04 (  -  ) &3.02e+04 (  -  )\\ 
\textit{roE100}    &3.17e+04   &3.56e+04 (\texttt{+}\textbf{12.3}\%)   &    -    &3.01e+04 (  -  ) &3.23e+04 (  -  )\\ 
\textit{roA100}    &3.44e+04   &3.32e+04 (- 3.5\%)   &    -    &3.32e+04 (  -  ) &3.29e+04 (  -  )\\ 
    \bottomrule
    \end{tabular}
\end{minipage}%
\end{table*}

\subsection{Heuristic Initialization of an Exact Algorithm}

The TSP-PC problem formulation defined by Eq. (\ref{eq:ProbForm}) is associated with binary variables $x_{ij}$ that defines which  edges $(i,j)\in \mathcal{E}$ to constitute the tour, and continuous variables $y_{im}$ related to payload evolution.
The quadratic term $x_{ij}y_{jm}$ in Eq. (\ref{eq:cargoEvolution}) sets up a Mixed Integer Nonlinear Programming (MINLP) formulation which is challenging to solve, especially considering the $\mathcal{NP}$-hard nature of the TSP-PC.
To address this, each $x_{ij}y_{im}$ is linearized to reduce problem complexity using the big $M$ method by introducing auxiliary variables $\lambda_{ijm} = x_{ij}y_{im}$ and four additional constraints:
\begin{subequations}  \label{eq:bigM}
\allowdisplaybreaks
\begin{align}  
    \label{eq:M1}   \lambda_{ijm}&  \leq y_{im}+M(1-x_{ij})     & \forall ~ i,j\in\overline{\mathcal{V}}, m\in\mathcal{H}\\
    \label{eq:M2}   \lambda_{ijm}&  \geq y_{im}-M(1-x_{ij})     & \forall ~  i,j\in\overline{\mathcal{V}}, m\in\mathcal{H}\\
    \label{eq:M3}   \lambda_{ijm}&  \leq Mx_{ij}                & \forall ~ i,j\in\overline{\mathcal{V}}, m\in\mathcal{H}\\
    \label{eq:M4}   \lambda_{ijm}&  \geq -Mx_{ij}               &  \forall ~  i,j\in\overline{\mathcal{V}}, m\in\mathcal{H}
\end{align} 
\end{subequations}

The value of $M$ in Eq. (\ref{eq:bigM}) must be chosen to be a sufficiently large positive number so that the auxiliary variable and additional constraints linearize $x_{ij}y_{im}$.
This is satisfied for $M > 2$ since the maximum payload of $m\in\mathcal{H}$ is at most 2 in the problem being studied.
By replacing $x_{ij}y_{im}$ in Eq. (\ref{eq:ProbForm}) with $\lambda_{ijm}$ and adding the constraints of Eq. (\ref{eq:bigM}), the MINLP is converted to a Mixed Integer Linear Programming problem (MILP) which can be solved effectively using commercially available solvers.

The MILP reformulation of Eq. (\ref{eq:ProbForm}) is modeled in a Matlab R2023b environment  and the \textit{intlinprog} solver of Gurobi 9.5.2 \citep{gurobi} is utilized on an AMD Ryzen 5600X CPU clocked at 3.7 GHz, paired with 128GB RAM.
The time limit of the exact solver is set to 5 minutes and the results are tabulated in Tables \ref{table1} and \ref{table2}.
The fourth column provides the solution found when running the exact algorithm \textit{standalone} without any heuristic initialization, where it is seen that for instances with 99 or more points, no feasible solution is found within the allocated computation time.
The fifth and sixth columns show the solution found when the exact solver uses the NN solution and ACHCI solution respectively for a warm start, along with the percentage reduction in solution cost compared with the standalone solution.

\subsection{Discussion}
The spatial characteristics of precedence constraints significantly affect the performance of the ACHCI algorithm since child nodes can only be inserted in feasible segments of the tour.
When parent nodes are located closer to the periphery of the point cloud, insertions onto the resulting convex hull are restricted to the region enclosed by the tour.
This avoids large cost increments when building the tour, and produces tours with low solution costs.
However, when parent nodes are centrally located, the insertions of child nodes from the periphery can result in large excursions from the tour, depending on the distance of the child node from its feasible insertion segment of the subtour.
When only a few child nodes remain to be inserted but their insertable segments happen to be far from them, the ACHCI solution cost is significantly affected.
In these cases, the NN heuristic performs better than the ACHCI.

\begin{figure}[b]
    \centering
    \vspace{-2mm}
        \includegraphics[trim =0mm 0mm 0mm 0mm, clip, width=0.8\linewidth]{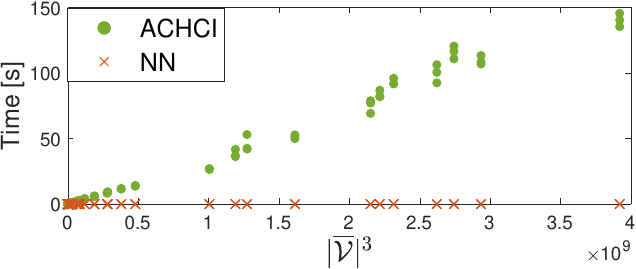}
  \caption{Computation time for the two heuristics}
  \label{4: PC-time} 
\end{figure}
A worst case complexity of $O(n^3)$ characterizes the ACHCI heuristic, as seen in Fig. \ref{4: PC-time} where the x axis is the number of TSP-PC locations cubed.
Because of the greedy nature of the NN heuristic, it almost instantaneously provides tours regardless of the number of points.
Considering the negligible computation time, it is worthwhile to compute the NN tour in addition to the ACHCI tour, and simply choose the tour with minimum cost between the two.

The initiation of the candidate subtour with the convex hull of the points is advantageous when precedence constraints do not exist, because points on the boundary of the convex hull are visited in the same cyclic order as their order in the optimal tour \citep{deineko1994convex}.
The ACHCI heuristic makes use of this property and therefore, the resulting solution acts as a good initialization point for the exact solver.
This is seen in Tables \ref{table1} and \ref{table2}, where the convergence of the exact solver is improved in every instance where the ACHCI tour is used as an initialization.
In most NN warm start cases as well, the convergence of the exact solver was improved.

\section{Conclusion}

The proposed ACHCI algorithm is a tour construction heuristic that accounts for precedence constraints by first drawing the initiating convex hull  over non-delivery locations and permitting subsequent insertions only in feasible segments that respect precedence constraints. 
Using the ACHCI heuristic solution to initialize an exact solver improved convergence in all cases, while the NN heuristic solutions did not always improve convergence.
The ACHCI heuristic solution outperforms the NN heuristic in all cases where delivery locations are centrally located, though it did not perform as well otherwise.
It is therefore highly applicable to spatial configurations of operations where pick-up locations are situated around the periphery to facilitate efficient collection and delivery locations are centrally situated for sorting, processing, or distribution.
Examples include layouts in flexible manufacturing systems or warehousing facilities where inventory shelves are at the perimeters of the facility, or in online retail and digital commerce settings where sorting facilities are located on the outskirts of a city and customers are within city premises.

\bibliography{ifacconf}             

\end{document}